\documentstyle[11pt,amstex]{article}

\newcommand{\CPtu}{{\Bbb C}P^2}
\newcommand{\cPtu}{{\Bbb C}P^2}

\newcommand{\lb}{\linebreak[0]}
\newcommand{\nline}{\ \\}

\newcommand{\pf}{{\bf Proof. }}

\newcommand{\spf}{{\bf Sketch of the proof. }}
\newcommand{\spfof}{{\bf Sketch of the proof of }}
\newtheorem{thm}{Theorem}[section]
\newtheorem{lem}[thm]{Lemma}
\newtheorem{rem}[thm]{Remark}
\newtheorem{prop}[thm]{Proposition}
\newtheorem{dfn}[thm]{Definition}
\newtheorem{cor}[thm]{Corollary}

\newcommand{\Qed}{\hfill \rule{.75em}{.75em}}

\newcommand{\C}{\:{\rlap {\raise 0.37ex \hbox{$\scriptstyle |$}}
\hskip -0.2em \hbox{\rm{C}}}}
\newcommand{\R}{\:{\rlap {\raise 0.37ex \hbox{$\scriptstyle |$}}
\hskip -0.2em \hbox{\rm{R}}}}
\newcommand{\Q}{\:{\rlap {\raise 0.37ex \hbox{$\scriptstyle |$}}
\hskip -0.2em \hbox{\rm{Q}}}}
\newcommand{\Z}{\:{\rlap {\raise 0.37ex \hbox{$\scriptstyle |$}}
\hskip -0.2em \hbox{\rm{Z}}}}

\newfont{\dbl}{msbm10 scaled\magstep1}

\newcommand{\bbar}{\overline}
\newcommand{\omgbar}{\overline{\omega}}
\newcommand{\Omgbar}{\overline{\Omega}}
\newcommand{\Jbar}{\overline{J}}
\newcommand{\sigstd}{\sigma_{std}}
\newcommand{\omgstd}{\omega_{std}}
\newcommand{\phit}{\varphi_{_t}}
\newcommand{\undx}{\underline{x}}
\newcommand{\dst}{\displaystyle}

\title{Connectedness of spaces of symplectic embeddings} 
\author{Paul Biran \\\\
	{\normalsize Preliminary version}}

\begin{document}

\maketitle

\begin{abstract}

We prove that the space of symplectic packings of $\CPtu$ by $k$ equal balls 
is connected for $3\leq k\leq 6$. The proof is based on Gromov-Witten 
invariants and on the inflation technique due to Lalonde and McDuff.

\end{abstract}

\section{Introduction}

Let $(M^{4},\omega)$ be a symplectic 4-manifold. Denote by 
$Emb(\coprod_{q=1}^{N}(B^{4}(\lambda_q),\omega_{std});\lb M,\omega)$
the space of symplectic embeddings of the disjoint union 
$\coprod_{q=1}^{N}B^{4}(\lambda_q)$ into $(M,\omega)$, where 
$B^{4}(\lambda)$ denotes the standard $4$-dimensional closed ball of radius
$\lambda$, equipped with the standard symplectic form 
$\omega_{std}=dx_1 \wedge dy_1 + dx_2 \wedge dy_2$. We endow 
these spaces with the $C^{\infty}$ topology.

In this paper we discuss $\pi_0(Emb)$, that is the set of 
symplectic isotopy classes of symplectic embeddings of balls into $M$. 

We shall be mainly concerned with the symplectic manifolds 
$(Int\,B^4(1),\lb \omgstd)$ and $(\CPtu,\sigstd)$, where $\sigma_{std}$ 
denotes the standard K\"{a}hler form on $\CPtu$, 
normalized such that $\int_{_{{\Bbb C}P^1}} \sigma_{std} = \pi$.

Quite a lot of research has been made concerning symplectic embeddings of 
balls into $(Int\,B^{4}(1),\omega_{std})$ and $(\CPtu,\sigma_{std})$.
The main directions of research are summarized bellow :
\begin{itemize}
	\item In \cite{M-P} McDuff and Polterovich found necessary and 
	sufficient conditions on the radii of the balls $\lambda_q$, for
	$\coprod_{q=1}^{N}B^4(\lambda_q)$ to embed symplectically 
	into ${\Bbb C}P^2$ or $Int\,B^4(1)$, when $1 \leq N \leq 9$, and for 
	$N=k^2,\;k\in{\Bbb N}$.

	\item McDuff proved in \cite{McD-Unique} that the spaces 
	$Emb(B^4(\lambda_1)\coprod B^4(\lambda_2); {\Bbb C}P^2)$ are connected 
	when non-empty.
\end{itemize}

In what follows we prove connectedness of the spaces of symplectic embeddings 
of 3,4,5 and 6 balls into $\CPtu$ and $Int\,B^4(1)$. 
As corollaries we prove some results on extensions of 
symplectic embeddings of balls to embeddings of larger balls.

The methods we use are based on Gromov's theory of pseudo-holomorphic curves, 
and following Lalonde and McDuff, on the inflation technique for correcting 
pseudo-isotopies to isotopies.

\nline
{\bf A remark on the notations.}\ \ \ Since our main objects of interest are
$(\CPtu,\sigma_{std})$ and $(Int\,B^4(1),\omega_{std})$,
$Emb(\coprod_{q=1}^{N}(B^{4}(\lambda_q),\omega_{std});M,\omega)$ will be 
denoted sometimes by $Emb(\lambda_1, \ldots ,\lambda_q)$, when no 
confusion may occur.

\section{Statement of the main results}
\label{sect-statement}

\subsection{Connectedness of spaces of symplectic embeddings}

\begin{prop}
\label{prop-pi0}

\nline
1) ${\dst Emb(\coprod_{1}^{3}B^4(\lambda);{\Bbb C}P^2,\sigma_{std})}$ is 
connected when it is non-empty (that is, when $\lambda<\frac{1}{\sqrt{2}}$).\\
2) ${\dst Emb(\coprod_{1}^{4}B^4(\lambda);{\Bbb C}P^2,\sigma_{std})}$ is 
connected when it is non-empty (that is, when $\lambda<\frac{1}{\sqrt{2}}$).\\
3) ${\dst Emb(\coprod_{1}^{5}B^4(\lambda);{\Bbb C}P^2,\sigma_{std})}$ is 
connected when it is non-empty (that is, when $\lambda<\sqrt{\frac{2}{5}}$).\\
4) ${\dst Emb(\coprod_{1}^{6}B^4(\lambda);{\Bbb C}P^2,\sigma_{std})}$ is 
connected when it is non-empty (that is, when $\lambda<\sqrt{\frac{2}{5}}$).\\
\end{prop}

The proof of this proposition appears in section \ref{sect-proofs-1}.
See also propositions \ref{prop-pi0-gen1} - \ref{prop-pi0-gen5} for slightly 
sharper results concerning balls of different radii, and arbitrary number 
of balls.\\

The next proposition gives some information on the connectedness of the space
of symplectic embeddings of arbitrary number of balls, provided that they are
small enough (compare with \cite{McD-Unique} section 3.6).

\begin{prop}
\label{prop-N-balls}  

${\dst Emb(\coprod_{1}^{N}B^4(\lambda_q);{\Bbb C}P^2,\sigma_{std})}$ 
is non-empty and connected when ${\dst \sum_{q=1}^{N} \lambda_q^2 < 1}$.
\end{prop}
Note that this proposition covers all the cases of embedding spaces 
of one and of two balls.

\subsection{Extensions of symplectic embeddings}

As a corollary of the previous propositions we prove:

\begin{cor}
\label{cor-extend}

1) Let $1\leq N\leq 6$, and let $\mu>0$ be such that \\
${\dst {Emb(\coprod_{1}^N B^4(\mu); \cPtu) \neq \emptyset}}$. 
Let $\lambda<\mu$, then every symplectic embedding 
${\dst {\coprod_{1}^N B^4(\lambda)\rightarrow (\cPtu,\sigstd)}}$, 
can be extended to a symplectic embedding \\  
${\dst {\coprod_{1}^N B^4(\mu)\rightarrow (\cPtu,\sigstd)}}$.\\
2) Let $1\leq N\leq 6$, and let $\mu>0$ be such that 
${\dst {Emb(\coprod_{1}^N B^4(\mu);\cPtu)\neq \emptyset}}$. Let $r<N$, then
every symplectic embedding 
${\dst {\coprod_{1}^r B^4(\mu)\rightarrow (\cPtu,\sigstd)}}$, 
can be extended to a symplectic embedding 
${\dst {\coprod_{1}^N B^4(\mu)\rightarrow (\cPtu,\sigstd)}}$.\\

\end{cor} 

The proof can be found in section \ref{sect-proofs-1}.

\section{Symplectic embeddings of balls and blow-ups}
\label{sect-blow-up}

In order to fix notations, we briefly review in this section
some results concerned with the relations 
between symplectic embeddings and blow-ups. The first results in this 
direction were established by McDuff in \cite{McD-Blow}, 
however, for our purposes it is more convenient to use the formulation 
appearing in \cite{M-P} and \cite{M-S-Intro}.

\subsection{Blowing-up in the symplectic category}

We start with the following definition 

\begin{dfn}

Let $(M,\omega)$ be a symplectic manifold, and let be $J$ an $\omega$-tamed 
almost complex structure on $M$. We say that $\omega$ is $J$-standard near 
$p \in M$ iff pair $(\omega,J)$ is diffeomorphic to the standard pair 
$(\omgstd,i)$ of ${\Bbb R}^{2n}$, near $p$.
\end{dfn}

Let $(M^{2n},\omega)$ be a closed symplectic manifold. 
Let $x\in M$, and choose an $\omega$-tamed almost complex structure $J$, 
for which $\omega$ is $J$-standard near $x$. Denote by 
$\Theta:(\bbar{M},\bbar{J}) \rightarrow (M,J)$ the complex 
blow-up of $M$ at x.
Let $\psi:(B^{2n}(\lambda) \rightarrow (M,\omega)$ be a symplectic embedding 
which is $(i,J)$ holomorphic near $0\in B^{2n}(\lambda)$, and such that 
$\psi(0)=x$.
The symplectic embedding $\psi$ gives rise to a symplectic form 
$\bbar{\omega}_{\psi}$ on $\bbar{M}$, in the cohomology class 
$[\bbar{\omega}_{\psi}]=[\Theta^*{\omega}]-\pi\lambda^2e$, where $e$ is the 
Poincar\'{e} dual to the homology class of the  exceptional divisor 
$\Sigma=\Theta^{-1}(x)$.
The form $\bbar{\omega}_{\psi}$ is called the symplectic blow-up of $\omega$ 
with respect to $\psi$, and is uniquely defined up to isotopy of forms.
For the construction of the form $\bbar{\omega}_{\psi}$, see 
\cite{M-P} section 5, or \cite{M-S-Intro} pages 223-225.

The same construction extends for arbitrary symplectic embeddings 
$\psi:B^{2n}(\lambda)\rightarrow(M,\omega)$
(that is, without assuming $\psi$ to be holomorphic near $0$), 
in view of the following lemma (see \cite{M-S-Intro}).

\begin{lem}
\label{lem-prtb-holo}

{\em (McDuff-Polterovich)} \ Let 
$\psi:(B^{2n}(\lambda),\omega_{std})\rightarrow(M,\omega)$ be a symplectic 
embedding, and let $J$ be an $\omega$-tamed almost complex structure on $M$,
such that $\omega$ is $J$-standard near x. Then for every compact subset 
$K\subseteq M\setminus\{x\}$ there exist a number $\delta \in (0,\lambda)$, a 
symplectic form $\omega'$ on $M$ which is isotopic to $\omega$ and is 
$J$-standard near $x$, and a symplectic embedding 
$\psi':(B^{2n}(\lambda),\omega_{std})\rightarrow(M,\omega')$ with the following
properties:\\
1) $\psi'\mid_{B(\delta)}$ is holomorphic.\\
2) $\psi'=\psi$ near $\partial B(\lambda)$.\\
3) $\omega'$ tames $J$ and coincides with $\omega$ on $K$.

\end{lem}

Thus, $\psi$ can always be slightly perturbed to be holomorphic near $0$, 
which enables us to define the blown-up form $\bbar{\omega}_\psi$ for 
arbitrary symplectic embeddings.

The next lemma describes to which extent does the definition of 
$\bbar{\omega}_\psi$ depend on the choices made 
(see \cite{M-S-Intro} lemma 6.47).

\begin{lem}

1) The isotopy class of the form $\bbar{\omega}_\psi$ depends only on 
the  embedding $\psi$ and on the germ of $J$ at $x$.\\
2) If two symplectic embeddings $\psi_1$ and $\psi_2$ are isotopic through a 
family of symplectic embeddings of $(B^{2n}(\lambda),\omega_{std})$ which take 
$0$ to $x$ then the corresponding blown-up forms are isotopic.

\end{lem}

All the above remains true for a family of balls as well, that is, every 
symplectic embedding 
$\psi=\coprod\psi_q:\coprod_{q=1}^{N} B^{2n}(\lambda_q)\rightarrow(M,\omega)$ 
gives rise to a 
symplectic form $\bbar{\omega}_{\psi}$ on $\bbar{M}$ (where $\bbar{M}$ is the 
blow-up of $M$ at the $N$ points $\psi_1(0),\ldots,\psi_N(0)$), 
and the previous theorems still hold in this case.

The next lemma shows that $\bbar{\omega}_{\psi}$ is pseudo-isotopic to a
form taming $\bbar{J}$ (see~\cite{M-S-Intro} proposition 6.46). 
More precisely

\begin{lem}
\label{lem-prep-1}
{\em (McDuff-Polterovich)} \ 
Let 
${\dst \psi:\coprod_{q=1}^N(B^{2n}(\lambda_q),\omega_{std})
\rightarrow(M,\omega)}$ be a symplectic embedding. 
Choose an $\omega$-tamed almost complex structure $J$ on $M$, for 
which $\omega$ is $J$-standard near $\psi_1(0),\ldots,\psi_N(0)$, and let 
$(\bbar{M},\bbar{J})$ be the complex blow-up at these
points. Then there exists a smooth family $\{\omgbar_t\}_{t\in[0,1]}$ 
of symplectic forms on $\bbar{M}$, such that $\bbar{\omega}_0$ tames 
$\bbar{J}$, and $\bbar{\omega}_1=\bbar{\omega}_{\psi}$.

\end{lem}

\nline
In what follows we'll be using the following key lemmas.
The first one is an obvious modification of lemma 6.48(ii) appearing in 
\cite{M-S-Intro}.

\begin{dfn}

Let $x_1,\ldots,x_N$ be distinct points in $M$, and 
$\undx=(x_1,\ldots,x_N)$. We say that a symplectic embedding 
$\psi=\coprod\psi_q:\coprod_{q=1}^{N}(B^{2n}(\lambda_q),\omega_{std})
\rightarrow(M,\omega)$ is $\undx$-normalized iff $\psi_q(0)=x_q$ 
for all q. Denote by 
$Emb_{\undx}(\coprod_{q=1}^{N}B^{2n}(\lambda_q);\lb M,\omega)$ 
the space of $\undx$-normalized symplectic embeddings of $N$ disjoint 
$2n$-balls of radii $\lambda_1,\ldots,\lambda_N$, into $M$.

\end{dfn} 

Denote by ${\cal C}_N(M)$ the configuration space of $M$, that is
$${\cal C}_N(M) := \{(x_1,\ldots,x_n) \mid 
x_q\in M,\; x_i\neq x_j \; for\; all\; i\neq j\}$$ 

\begin{lem}
\label{lem-normal-emb}

If ${\cal C}_N(M)$ is simply connected then the inclusion 
$$Emb_{\underline{x}}(\coprod_{q=1}^{N}B^{2n}(\lambda_q);M,\omega)
\rightarrow Emb(\coprod_{q=1}^{N}B^{2n}(\lambda_q);M,\omega)$$ 
induces an isomorphism on the corresponding spaces of path components.

\end{lem}
Note that for $M=\CPtu$ or $Int\,B^4(1)$ for example, ${\cal C}_N(M)$ is 
simply connected for all~$N$, hence the above lemma implies that if two 
symplectic embeddings 
$\varphi,\psi \in Emb_{\undx}(\coprod_{q=1}^N B^{4}(\lambda_q); 
\CPtu,\sigstd)$ are symplectically isotopic then they are also 
isotopic through an $\undx$-normalized family of symplectic embeddings.

\begin{lem}
\label{lem-isotopy}

Let ${\dst \phi=\coprod_{q=1}^{N}\phi_q}$ and 
${\dst \psi=\coprod_{q=1}^{N}\psi_q}$ be two $\undx$-normalized symplectic 
embeddings ${\dst \coprod_{q=1}^{N}(B^{4}(\lambda_q),\omega_{std})
\rightarrow(\cPtu,\sigstd)}$. Then $\phi$ and $\psi$ are symplectically 
isotopic iff the blown-up forms $\bbar{\omega}_{\phi}$ and 
$\bbar{\omega}_{\psi}$ are isotopic.

\end{lem}
This is a straightforward consequence of the previous lemma and of 
propositions 2.6 and 2.7 from \cite{McD-Unique}.

\nline

The following proposition shows that from our point of view there is no
essential difference between $M={\Bbb C}P^2$ and $M=Int\,B^4(1)$,
(see~\cite{McD-Unique} proposition 2.8).

\begin{prop}
\label{prop-ball-cp2}

{\em (McDuff)} \ Denote by 
$\iota:Int B^4(1)\rightarrow{\Bbb C}P^2\setminus{\Bbb C}P^1$ the map 

$$\iota(z,w)=[\sqrt{1-\mid z\mid^2-\mid w \mid^2}:z:w] \;\;\;\;\;\;  
(z,w)\in Int\,B^4(1)\subseteq{\Bbb C}\times{\Bbb C}$$
Then\\
1) $\iota$ is a symplectomorphism 
$(Int B^4(1),\omega_{std})\rightarrow
({\Bbb C}P^2\setminus{\Bbb C}P^1,\sigma_{std})$.\\
2) $\iota$ induces an isomorphism between 
$\pi_{0}(Emb(\coprod B^4(\lambda_q); Int\,B^4(1),\omega_{std}))$ and \\
$\pi_{0}(Emb(\coprod B^4(\lambda_q);{\Bbb C}P^2,\sigma_{std}))$.

\end{prop}

In view of this proposition one can replace $({\Bbb C}P^2,\sigma_{std})$ by 
$(Int\,B^4(1),\omega_{std})$ in every place in the following,
unless otherwise stated.

\section{Pseudo-holomorphic spheres}

In this section we use pseudo-holomorphic theory to prove two key propositions 
which will be used later to prove the connectedness of certain symplectic 
embedding spaces.

\begin{prop}
\label{prop-holo-sphere}

Let $(M^4,\omega_0)$ be a symplectic closed 4-manifold. Let $S_0\subseteq M$ be
a symplectically embedded 2-sphere which represents the 2-homology class 
$A=[S_0]\in H_2(M;{\Bbb Z})$, with $A\cdot A=0$.
Suppose that there does not exist any $B\in H_2(M;{\Bbb Z})$ such that $A=2B$.
Let $J_0$ be an $\omega_0$-tamed almost complex structure on M for which 
$S_0$ is $J_0$-holomorphic and suppose that $J_0$ is positive (that is, 
there do not exist any $J_0$-holomorphic spheres with non-positive 1'st 
Chern number. See remark bellow). Let $\{\omega_t\}_{t\in[0,1]}$ be a pseudo-isotopy of 
symplectic forms on M, starting at $\omega_0$. Then there exists a path of 
$\omega_t$-tamed almost complex structures $J_t$ starting at $J_0$, and a 
smooth family of embedded $J_t$-holomorphic spheres $S_t$, starting at $S_0$, 
and representing the homology class $A$.

\end{prop}

\begin{rem}
\label{rem-holo-sphere}

\nline
1) The condition on the positiveness of $J_0$ can weakened by requiring that 
there do not exist any $J_0$-holomorphic $A$-cusp curves (with components of 
genus 0), having some component with non-positive 1'st Chern number.\\
2) Note that $J_0$ can always be slightly perturbed in such a way that it 
becomes positive, since the set of positive almost complex structures on $M^4$ 
is dense in the space of all almost complex structures tamed by $\omega_0$.
Note also that given $K>0$, the set of all $\omega$-tamed 
almost complex structures $J$ on $M^4$ satisfying the following property:
\begin{quote}
for each $B\in H_2(M;{\Bbb Z})$ with $c_1(B)\leq 0$ and $\omega(B)\leq K$, 
there are no $J$-holomorphic $B$-spheres
\end{quote}
is open dense and path connected.
In particular, the set of all $\omega$-tamed almost complex structures $J$, 
for which there do not exist any $J$-holomorphic cusp $A$-curves having a 
component with non-positive 1'st Chern number, contains an open dense and 
path connected subset. 

\end{rem}

\spfof {\bf proposition \ref{prop-holo-sphere}.\ \ }
Consider the Gromov-Witten invariant $\Phi_{A}([pt])$.
The idea is to show that $\Phi_{A,J_0}([pt])=1$, then to take a generic path 
$J_t$ of $\omega_t$-tamed almost complex structures starting at $J_0$, and 
use the independence of $\Phi_A$ under deformations to show that 
$\Phi_{A,J_t}([pt])=\Phi_{A,J_0}([pt])=1$. This will prove that there exists 
a $J_t$-holomorphic sphere in the class $A$ for all $t$. 
These spheres are clearly embedded by the adjunction formula.
The homological condition that $A$ is not divisible by 2 is imposed 
in order to control the behavior of $A$-cusp curves, thus insuring the 
possibility of using Gromov-Witten invariants.

Finally we have to prove the existence of a smooth family $S_t$ of such 
spheres. To prove this, note that since $\Phi_{A,J_t}([pt])=1$, there exists 
exactly one $J_t$-holomorphic $A$-sphere through every point which does not 
lie on a $J_t$-holomorphic $A$-cusp curve. Since the images of the 
$J_t$-holomorphic $A$-cusp curves have at least codimension 2 in $M$, we can
choose a smooth path $p_t$ of points in $M$ starting at a point on $S_0$, 
such that for all $t$ there are no $J_t$-holomorphic $A$-cusp curves 
passing through $p_t$. This implies that for
all t there exists exactly one $J_t$-holomorphic $A$-sphere $S_t$ 
through $p_t$.
\Qed

\subsection{Inflation}
Following Lalonde and McDuff (see \cite{L-M-J}, \cite{McD-Notes} \ lemma 3.7), 
we use the so called inflation procedure to construct $J$-semi positive 
closed 2-forms in suitable cohomology classes. 
More precisely, given a symplectically embedded 2-sphere, 
with self intersection number 0, and an almost complex structure $J$, for 
which this sphere is holomorphic, this procedure enables us to build a 
closed 2-form $\tau$, which is Poincar\'{e} dual to the class of the given 
sphere, and which {\em weakly tames} $J$, that is $\tau(X,JX) \geq 0$ for all
$X\in TM$.\\
We call the following, the {\em inflation lemma}.

\begin{prop}
\label{prop-inflation}

{\em (inflation lemma)} \ Let $(M^4,\omega)$ be a 
symplectic closed 4-manifold. Let $S\subseteq M$ be 
a symplectically embedded 2-sphere with $S\cdot S=0$. Then for every 
$\omega$-tamed almost complex structure $J$ for which $S$ is 
$J$-holomorphic, there exist a closed 2-form $\tau_{_J}$, weakly taming $J$,
supported in an arbitrarily small neighborhood of $S$, and representing the 
cohomology class Poincar\'{e} dual to $[S]\in H_2(M;{\Bbb Z})$. 
Moreover, if $\omega_t$ is a pseudo-isotopy of symplectic forms, $\{J_t\}_t$ 
a smooth path of $\omega_t$-tamed almost complex structures, and $S_t$ a 
smooth family of embedded 
$J_t$-holomorphic 2-spheres with $S_t\cdot S_t=0$, then $\tau_{_{J_t}}$ can be 
chosen to depend smoothly on t.

\end{prop}

\spf
>From our assumptions it follows that $S$ has trivial normal bundle.
The idea is to construct a tubular neighborhood $W(S)$ of $S$, consisting of 
2-discs in the normal direction to $S$, and of $A$-holomorphic spheres in the 
horizontal direction (that is, sections of the tubular neighborhood), where 
$A=[S]\in H_2(M;{\Bbb Z})$. This is done by standard arguments from 
pseudo-holomorphic theory.\\
The almost complex structure $J$ induces an orientation on these 2-discs,
since $S$ is $J$-holomorphic.
Choose one such oriented 2-disc $D$, which is transverse to the $A$-spheres, 
and consider the projection $\pi:W(S)\rightarrow D$.
Choose now a bump 2-form $\alpha$ on $D$, with small enough support, such 
that $\int_{_D} \alpha = 1$. The form $\tau_{_J}=\pi^*\alpha$ is what we need.
\Qed 

\nline\noindent
Combining propositions \ref{prop-holo-sphere} and \ref{prop-inflation} 
we have

\begin{prop}
\label{prop-holo-combin}

Let $(M^4,\omega_0)$ be a symplectic closed 4-manifold. Let $S_0\subseteq M$ 
be an $\omega_0$-symplectically embedded 2-sphere with $S_0\cdot S_0=0$.
Suppose that there does not exist any $B\in H_2(M;{\Bbb Z})$ such that 
$[S_0]=2B$. Let $J_0$ be a positive almost complex structure tamed by 
$\omega_0$ (see remark bellow), such that $S_0$ is $J_0$-holomorphic.
Let $\{\omega_t\}_{t\in[0,1]}$ be a pseudo-isotopy of symplectic forms 
starting at $\omega_0$. Then there exists a path of $\omega_t$-tamed almost 
complex structures $J_t$, and a smooth path of closed 2-forms $\tau_t$  
weakly taming $J_t$, such that $[\tau_t]=PD([S_0])$ for all $t$.
Moreover, the above remains true if we have a family 
$S_0^{1},\ldots,S_0^{r}$ of embedded $J_0$-holomorphic spheres with self
intersection numbers 0, that is, there 
exists a smooth path of closed 2-forms $\tau_t^{1},\ldots,\tau_t^{r}$  
weakly taming $J_t$, such that $[\tau_t^{i}]=PD([S_0^{i}])$ 
for all $1\leq i\leq r$.

\end{prop}

\begin{rem}
\label{rem-holo-combin}

The same remarks following proposition \ref{prop-holo-sphere}, apply here
as well. 

\end{rem}

\section{$\pi_0$ of spaces of symplectic embeddings}

In this section we describe the method which we use in order to prove 
the connectedness of certain symplectic embedding spaces.

Let $\varphi_{_{-1}}$, $\varphi_{{_1}}$ be two symplectic embeddings of 
balls into ${\Bbb C}P^2$. 
As was explained in section \ref{sect-blow-up}, they give rise to 
symplectic forms $\bbar{\Omega}_{-1}$ and $\bbar{\Omega}_{1}$, on a certain 
blow-up of ${\Bbb C}P^2$. 
By lemma \ref{lem-isotopy},
in order to show that $\varphi_{_{-1}}$ and $\varphi_{_1}$ are isotopic we 
need to show that $\bbar{\Omega}_{-1}$ and $\bbar{\Omega}_{1}$ are isotopic.

By lemma \ref{lem-prep-1}, each of $\bbar{\Omega}_{-1}$ and 
$\bbar{\Omega}_{1}$ is 
pseudo-isotopic to forms which tame a common almost complex structure, hence
$\bbar{\Omega}_{-1}$ and $\bbar{\Omega}_{1}$ are mutually pseudo-isotopic. 
Denote the path connecting them by $\{\bbar{\Omega}_t\}_{t\in[-1,1]}$.

The idea is to correct the pseudo-isotopy $\{\bbar{\Omega}_t\}_{t\in[-1,1]}$ to 
an isotopy connecting $\bbar{\Omega}_{-1}$ and $\bbar{\Omega}_{1}$ 
by taking a positive linear combination 
$x(t)\bbar{\Omega}_t + y(t)\tau_t$, 
of $\bbar{\Omega}_t$ and $\tau_t$, where $\tau_t$ is a suitable closed 2-form 
obtained by using the inflation lemma. We give now the precise details of the 
method.

\subsection{From pseudo-isotopy to isotopy}

We denote ${\Bbb C}P^2$ by $V$, to emphasize that the following method 
can work in some more general situations.

Fix $p_1,\ldots,p_N\in V$, points in general position, and set 
$\underline{p}=(p_1,\ldots,p_N)$. Let $\varphi_{_{-1}}$ , $\varphi_{_1}$ be 
two $\underline{p}$-normalized symplectic embeddings
$$\varphi_{_i}=
\coprod_{q=1}^{N}\varphi_{_{i,q}}:\coprod_{q=1}^{N}
(B^4(\lambda_q),\omega_{std})
\rightarrow ({\Bbb C}P^2,\sigma_{std}) \;\;\;\;\;  i=-1,1$$
Without loss of generality (see lemma \ref{lem-prtb-holo}) we can assume 
that $\varphi_{_{-1}},\varphi_{_1}$ are $(i,J_s)$-holomorphic near the 
centers of the balls, where $J_s$ is the standard complex structure on 
$V={\Bbb C}P^2$. 
Let $(\bbar{V},\bbar{J_s}) \stackrel{\Theta}{\rightarrow} (V,J_s)$ be the 
complex blow-up at $p_1,\ldots,p_N$. 
Fix a copy $\Sigma_L$ of ${\Bbb C}P^1 \subset {\Bbb C}P^2$ and let 
$\Sigma_{E_q}=\Theta^{-1}(p_q)$ \  $q=1,\ldots,N$ be the exceptional 
divisors. Let $L=[\Sigma_L]$, $E_q=[\Sigma_{E_q}]$ be the corresponding 
classes in 
$H_2(M;{\Bbb Z})$, and $l=PD(L),\; e_q=PD(E_q) \in H^2(M;{\Bbb Z})$ be their 
Poincar\'{e} duals. Finally, denote by $\Omgbar_{-1},\; \Omgbar_1$ the 
blow-up of the form $\sigma_{std}$ 
with respect to the symplectic embeddings $\varphi_{_{-1}},\; \varphi_{_1}$.

First we need a slightly more detailed version of lemma \ref{lem-prep-1}, 
which will be called the {\em preparation lemma}.
For technical reasons it is more convenient for us to replace the
blown-up forms $\bbar{\Omega}_i$, by $\frac{1}{\pi}\Omgbar_i$, so that 
their cohomology classes are $l-\sum_{q=1}^{N}\lambda_q^2e_q$, rather than 
$\pi(l-\sum_{q=1}^{N}\lambda_q^2e_q)$.\\\\

\begin{lem}
\label{lem-prep}
{\em (preparation lemma)}\ In the above situation, 
there exists a pseudo-isotopy 
$\{\Omgbar_t\}_{t\in[-1,1]}$ connecting $\Omgbar_{-1}$ and $\Omgbar_1$, 
such that: \\
1) $[\Omgbar_t]=l-\sum_{q=1}^{N}\alpha_q(t)^2 e_q$.\\
2) $\alpha_q(-1)=\alpha_q(1)=\lambda_q$.\\
3) $0<\alpha_q(t)\leq\lambda_q$ for all $t\in[-1,1]$.\\
4) $\Omgbar_t$ is symplectic on $\Sigma_{E_q}$ for all $q$.\\
5) $\Omgbar_0$ tames $\Jbar_s$.\\
Moreover, if $\lambda_i=\lambda_j$ for some $i\neq j$ then $\alpha_i,\alpha_j$
can be chosen to satisfy $\alpha_i(t)=\alpha_j(t)$ for all $t\in[-1,1]$.

\end{lem}
  
\spf
Since this lemma is almost the same as \ref{lem-prep-1}, we add here only the 
details which differ from \ref{lem-prep-1}.

According to lemma \ref{lem-prep-1} there exist pseudo-isotopies
$\{\bbar{\Omega}'_t\}_{t\in[-1,0]}$ and
$\{\bbar{\Omega}''_t\}_{t\in[0,1]}$ on $\bbar{V}$, such that 
$\bbar{\Omega}'_{-1}=\bbar{\Omega}_{-1}, 
\;\; \bbar{\Omega}''_1=\bbar{\Omega}_{1}$  
(where $\bbar{\Omega}_{-1},\;\bbar{\Omega}_{1}$ are the given blown-up forms) 
and such that
\begin{displaymath}
[\bbar{\Omega}'_t]=l-\sum_{q=1}^{N}\alpha'_q(t)^2 e_q \ \ \ \ \
[\bbar{\Omega}''_t]=l-\sum_{q=1}^{N}\alpha''_q(t)^2 e_q \ \ \ \ \ 
 \alpha'_q(-1)=\alpha''_q(1)=\lambda_q
\end{displaymath}

As one can see from the proof of lemma \ref{lem-prep-1} 
(appearing in \cite{M-P}, proposition 2.1.B), the forms $\bbar{\Omega}'_t$ and 
$\bbar{\Omega}''_t$ are obtained by shrinking the embedded balls to 
balls of radii $\alpha_q(t)$, and then blowing-up the form $\sigstd$ 
with respect to the shrinked 
embeddings. Since the shrinking of each of the balls can be made independently 
of the others we can control the values of $\alpha_q(t)$ in such a way that 
$0<\alpha'_q(t),\alpha''_q(t)\leq\lambda_q$, and even that
$\alpha'_i(t)=\alpha'_j(t)$,
$\alpha''_i(t)=\alpha''_j(t)$ when 
$\lambda_i=\lambda_j$.
Furthermore, the shrinking can be made in such a way that when $t$ 
reaches $0$, $\varphi_{_{-1}}$ and $\varphi_{_1}$ are holomorphic on 
$B^4(\alpha_q(t))$, hence the blown-up forms 
$\bbar{\Omega}'_0$ and $\bbar{\Omega}''_0$ both tame $\bbar{J}_s$.
This implies that $\bbar{\Omega}'_0$ and $\bbar{\Omega}''_0$ 
can be joined by a linear path of symplectic forms, namely 
$s\bbar{\Omega}'_0+(1-s)\bbar{\Omega}''_0$, $s\in[0,1]$.

We have now a path joining $\bbar{\Omega}'_{-1}$ with 
$\bbar{\Omega}''_1$, formed from 3 parts. Re-parameterizing it and
smoothing in the corners (that is, at $\bbar{\Omega}'_0$ and 
$\bbar{\Omega}''_0$) we get the desired pseudo-isotopy.

\Qed

\nline
Applying the above lemma gives us a pseudo-isotopy 
$\{\Omgbar_t\}_{t\in[-1,1]}$ connecting 
$\bbar{\Omega}_{-1}$ with $\bbar{\Omega}_1$, such that $\bbar{J}_s$ is tamed by 
$\bbar{\Omega}_0$. The next step is to correct this 
pseudo-isotopy to an isotopy.

Consider suitable (the meaning of ``suitable'' will become clear immediately)
$\bbar{J}_s$-holomorphic spheres $S^{1},\ldots,S^{r}$ in $\bbar{V}$ 
with self intersection numbers 0. 
Using proposition~\ref{prop-holo-combin} we obtain a 
path $\bbar{J}_t$, of almost complex structures tamed by $\bbar{\Omega}_t$, 
and a path of $\bbar{J}_t$-weakly taming closed 2-forms 
$\tau_t^{1},\ldots,\tau_t^{r}$, in the classes $s^{1},\ldots,s^{r}$, 
where $s^{j}$ is Poincar\'{e} dual to $[S^{j}]$.

The final step is crucial, and it reveals the meaning of the word 
``suitable''.\\
Suppose we have smooth families of real numbers, $x(t)>0$ and 
$y^{1}(t),\ldots,y^{r}(t) \lb \geq 0$, such that $x(-1)=x(1)=1,\;\;
 y^j(-1)=y^j(1)=0$ and such that

$$x(t)(l-\sum_{q=1}^{N}\alpha_q(t)^2 e_q)+\sum_{j=1}^{r}y^j(t)s^{j}=
l-\sum_{q=1}^N\lambda_q^2 e_q$$
If such families of numbers exist then consider the forms,

$$\bbar{\omega}_t:=x(t)\bbar{\Omega}_t+\sum_{j=1}^{r}y^{j}(t)\tau_t$$
Clearly $\bbar{\omega}_t$ is symplectic for all $t$, because $\tau_t$ weakly 
tames $\bbar{J}_t$. From our choice of $x(t),y^j(t)$, it follows that 
$\bbar{\omega}_t$ is an isotopy connecting 
$\bbar{\omega}_{-1}=\bbar{\Omega}_{-1}$ and 
$\bbar{\omega}_1=\bbar{\Omega}_1$, hence by lemma \ref{lem-isotopy}, the 
symplectic embeddings $\varphi_{_{-1}}$ and $\varphi_{_1}$ are 
isotopic through a family of symplectic embeddings.

\nline
In summary, the above method works when we can find 
$\bbar{J}_s$-holomorphic spheres 
$S^{1},\ldots,S^{r}$ in $\bbar{V}$ with self intersection numbers 0,
such that there exist smooth families of real numbers 
$x(t)>0 , \; y^j(t)\geq 0$ \ 
$j=1,\ldots,r$, satisfying the following cohomological condition:

$$x(t)[\bbar{\Omega}_t]+\sum_{j=1}^r y^j(t)PD([S^{j}])=
l-\sum_{q=1}^N \lambda_q^2 e_q $$
$$ x(-1)=x(1)=1,\;\; y^j(-1)=y^j(1)=0$$
Examples of applying this method for concrete cases appear in 
section \ref{sect-proofs-1}.

\section{Proofs of the main propositions}
\label{sect-proofs-1}

{\bf Proof of proposition \ref{prop-pi0}(2)}

Let $\varphi_{_{-1}},\varphi_{_1}:\coprod_{1}^{4}(B^4(\lambda),\omega_{std})
\rightarrow({\Bbb C}P^2,\sigma_{std})$ be two symplectic embeddings with 
$\varphi_{_{-1,_q}}(0)=\varphi_{_{1,_q}}(0) \; q=1,\ldots,4$, and without 
loss of generality assume that these 4 points are in general position.
Denote by $\CPtu_4$ the blow-up of $\CPtu$ at the 4 points 
$\varphi_{_{1,_1}}(0),\ldots,\varphi_{_{1,_4}}(0)$, and 
let $\bbar{\Omega}_{-1},\bbar{\Omega}_1$ be the blow-ups of $\sigstd$ with 
respect to $\varphi_{_{-1}},\varphi_{_1}$ on $\bbar{V}={\Bbb C}P^2_4$. 

By the preparation lemma (\ref{lem-prep}), there exists a 
pseudo-isotopy $\{\bbar{\Omega}_t\}_{t\in[-1,1]}$ connecting 
$\bbar{\Omega}_{-1}$ and $\bbar{\Omega}_1$, such that 
$[\bbar{\Omega}_t]=l-\alpha(t)^2\sum_{i=1}^4 e_i$, where 
$0<\alpha(t)\leq\lambda$ 
for all $t$, $\alpha(-1)=\alpha(1)=\lambda$, and such that $\Omgbar_0$ 
tames $\bbar{J}_s$, where $\Jbar_s$ is the lift of the standard complex 
structure of $\CPtu$ to $\CPtu_4$. Note that $\Jbar_s$ is positive.

Consider a $\bbar{J}_s$-holomorphic sphere $S_0$, in the homology 
class $A:=2L-\sum_{i=1}^4 E_i$ (that is, the proper transform of a conic 
through the 4 blown-up points). 
By proposition \ref{prop-holo-combin}, there exists 
a path $\{\Jbar_t\}_{t\in[-1,1]}$, of almost complex structures tamed by 
$\bbar{\Omega}_t$, and a smooth family of closed 2-forms $\tau_t$, weakly 
taming $\bbar{J}_t$, such that $[\tau_t]=2l-\sum_{i=1}^4 e_i$.

We find now smooth families of real numbers $x(t)>0$ and $y(t)\geq0$, such that
$x(-1)=x(1)=1,\;\; y(-1)=y(1)=0$, and such that 
$\bbar{\omega}_t:=x(t)\bbar{\Omega}_t+y(t)\tau_t$ lies in the cohomology class 
$l-\lambda^2\sum_{i=1}^4 e_i$, for all t. The forms $\bbar{\omega}_t$ are 
symplectic because 
$\tau_t$ weakly tames $\bbar{J}_t$.
Once we prove the existence of the coefficients $x(t), y(t)$ we are done, 
since $\bbar{\omega}_t$ is an isotopy connecting 
$\bbar{\omega}_{-1}=\bbar{\Omega}_{-1}$ with 
$\bbar{\omega}_1=\bbar{\Omega}_1$, hence by lemma \ref{lem-isotopy} 
$\varphi_{_{-1}}$ and $\varphi_{_1}$ are isotopic. The connectedness of 
$Emb(\coprod_{1}^4 B^4(\lambda); \CPtu,\sigma_{std})$ follows now from 
lemma~\ref{lem-normal-emb}, and the remarks after it.

It remains to prove the existence of the coefficients $x(t),y(t)$.
They must satisfy 

$$x(t)[{\Omgbar}_t]+y(t)[\tau_t]=l-\lambda^2\sum_{i=1}^4 e_i$$
This is equivalent to the system of linear equations

$$x(t)+2y(t)=1$$
$$\alpha(t)^2 x(t)+y(t)=\lambda^2$$
which has the solution 
$$x(t)=\frac{1-2\lambda^2}{1-2\alpha(t)^2} \;\;\;\;\;
y(t)=\frac{\lambda^2-\alpha(t)^2}{1-2\alpha(t)^2}$$
Indeed, $x(t)>0,\;\; y(t)\geq0$ and $x(-1)=x(1)=1,\;\; y(-1)=y(1)=0$.
\Qed

\nline
{\bf Proof of proposition \ref{prop-pi0}(3)}

Essentially it is the same as the proof of (2), but here we make inflation 
with 5 curves in the classes 
$$A_i:=2L-\sum_{q=1,q\neq i}^5 E_q, \;\;\;i=1,\ldots,5$$ 
to obtain closed 2-forms $\tau_t^{i}$ in the classes 
$$[\tau_t^{i}]=2l-\sum_{q=1,q\neq i}^5 e_q$$
Then we correct the pseudo-isotopy $\Omgbar_t$, where
$[\Omgbar_t]=l-\alpha(t)^2\sum_{q=1}^5 e_q$, by finding 
\mbox{$x(t)>0$}, \mbox{$y(t)\geq0$} such that 
$x(t)[\Omgbar_t]+y(t)\sum_{q=1}^5[\tau_t^{q}]=l-\lambda^2\sum_{i=1}^5 e_i$.
This is done by solving the system 
$$x(t)+10y(t)=1$$
$$\alpha(t)^2 x(t)+4y(t)=\lambda^2$$
which has the solution
$$x(t)=\frac{2-5\lambda^2}{2-5\alpha(t)^2}$$
$$y(t)=\frac{\lambda^2-\alpha(t)^2}{2(2-5\alpha(t)^2)}\;\;\; .$$
It is easy to see that $x(t)$, $y(t)$ satisfy the desired conditions.
\Qed

\nline
The proofs of proposition \ref{prop-pi0}(1),(4) require some technical details
so we give here only outlines of their proofs. 

\nline
{\bf Outline of the proof of proposition \ref{prop-pi0}(1)}
The idea is to make inflation with conics as we did for 4 balls.
The problem is that for making inflation we need curves with self 
intersection numbers zero, while here having only 3 points blown-up, 
any two conics through them intersect. 
To overcome this, we artificially blow-up another point and 
make inflation with conics through the 4 blown-up points. After correcting 
the pseudo-isotopy we blow it down at the 4'th exceptional divisor to achieve 
the isotopy we need. We give now more precise arguments, realizing 
this idea.\\

Consider two symplectic embeddings 
$\varphi_{_{-1}},\varphi_{_1}:\coprod_1^3 B^4(\lambda)\rightarrow \CPtu$, 
with $\varphi_{_{-1_q}}(0)=\varphi_{_{1_q}}(0)$, and let
$\Omgbar_{-1},\Omgbar_1$ be the corresponding blown-up form on 
$\bbar{V}=\CPtu_3$.
As before, we get a pseudo-isotopy $\{\Omgbar_t\}_t$ in the cohomology class
$[\Omgbar_t]=l-\alpha(t)^2\sum_{q=1}^3e_q$, connecting $\Omgbar_{-1}$ with 
$\Omgbar_1$.

Choose a point $q\in V$ which does not lie in the image of $\varphi_{_{-1}}$.
Clearly $\varphi_{_1}$ can be isotoped in such a way that $q$ does not lie in 
the image of $\varphi_{_1}$ too. Consider now the point 
$p=\Theta^{-1}(q) \,\in \bbar{V}$. It follows from the proof of the preparation
lemma that we can assume that $\Omgbar_t$ is standard near $p$ for all $t$, 
that is $\Omgbar_t=\Theta^*\sigma_{std}$ near p.
Choose now a very small $\beta>0$ for which there exist an 
$(i,\bbar{J}_s)$-holomorphic embedding\\
$\phi:B^4(1+\epsilon)\rightarrow\bbar{V}$, normalized at $p$ and  
such that $\phi^*\Omgbar_t=\beta^2\omega_{std}$.

Next, blow-up the forms $\Omgbar_t$ with respect to $\phi$ to obtain 
a new pseudo-isotopy of symplectic forms $\bbar{\Omgbar}_t$ in the cohomology 
class $[\bbar{\Omgbar}_t]=l-\alpha(t)^2\sum_{q=1}^3 e_q \lb - \beta^2 e_4$.

We now correct $\bbar{\Omgbar}_t$ as we did in the proof on 4 balls by making 
inflation with conics through the 4 blown-up points to produce the 
closed 2-forms $\tau_t$ in the cohomology class  
$[\tau_t]=2l-\sum_{q=1}^4e_q$, and then take 
$\bbar{\omgbar}_t:=x(t)\bbar{\Omgbar}_t+y(t)\tau_t$, where 
$$x(t)+2y(t)=1$$ $$\alpha(t)^2x(t)+y(t)=\lambda^2$$
Note that $\bbar{\omgbar}_t$ is not an isotopy yet, since 
$[\bbar{\omgbar}_t]=l-\lambda^2\sum_{q=1}^3e_q-\mu(t)e_4$, where 
$\mu(t)=x(t)\beta^2+y(t)$. The point is that after blowing-down 
$\bbar{\omgbar}_t$ back to $\CPtu_3$, it becomes an 
isotopy $\omgbar_t$, connecting $\Omgbar_{-1}$ with $\Omgbar_1$, hence 
the two embeddings $\varphi_{_{-1}}$ and $\varphi_{_1}$ are symplectically 
isotopic.
\Qed

\nline
{\bf Outline of the proof of proposition \ref{prop-pi0}(4)}

The proof is very similar to the proof of \ref{prop-pi0}(1), only that here 
we make inflation with quintics on $\CPtu_7=(\CPtu_6)_{_1}$ rather then 
with conics. The relevant homology class here is $A=5L-2\sum_{q=1}^6E_q-E_7$.
\Qed

\nline\nline
{\bf Proof of proposition \ref{prop-N-balls}}

Here we make inflation with curves in the classes $A_i:=L-E_i$,\ \ 
$i=1,\ldots,N$ on the blow-up of $\CPtu$ at $N$ points, $\bbar{V}=\CPtu_N$.
The standard complex structure $\Jbar_s$ on $\CPtu_N$ is no longer positive 
when \mbox{$N\geq9$}, however following the remark after proposition 
\ref{prop-holo-sphere}, one shows that there do not exist any 
$\Jbar_s$-holomorphic $A_i$-cusp curves which have a component with 
non-positive 1'st Chern number. 

The proof now proceeds as in \ref{prop-pi0}(2). The correction of the 
pseudo-isotopy is made by finding $x(t),y_1(t),\ldots,y_N(t)$ such that

$$x(t)(l-\sum_{q=1}^N\alpha_q(t)^2 e_q)+\sum_{q=1}^N y_q(t)(l-e_q)=
l-\sum_{q=1}^N\lambda_q^2 e_q$$
One easily computes that 
$$x(t)=\frac{1-\sum_{q=1}^N\lambda_q^2}{1-\sum_{q=1}^N \alpha_q(t)^2}, 
\;\;\;\;\;y_j(t)=\lambda_j^2-\alpha_j(t)^2x(t)$$
are suitable for our purposes.\\

The fact that $Emb(\coprod_{q=1}^N B^4(\lambda_q);\CPtu,\sigstd)$ is non-empty 
follows from Nakai-Moishezon criterion, by which it is easy to show that the 
cohomology class $l-\sum_{q=1}^N \lambda_q^2 e_q$ is  K\"{a}hler when 
$\sum_{q=1}^N \lambda_q^2<1$ (see \cite{M-P} theorem 1.3.D).

\Qed

\nline\nline
{\bf Proof of corollary \ref{cor-extend}}

The proof is based on the following proposition
\begin{prop}
\label{prop-isotop-extend}

Suppose that ${\dst Emb(\coprod_{q=1}^N B^4(\mu_q);\cPtu,\sigma_{std})}$ is 
non-empty and connected.\\
1) Let $0<\lambda_i\leq\mu_i$ for all $i$, then 
${\dst Emb(\coprod_{q=1}^N B^4(\lambda_q);\cPtu,\sigma_{std})}$ is 
connected iff the restriction map
$$Emb(\coprod_{q=1}^N B^4(\mu_q);\cPtu,\sigma_{std})\stackrel{Res}
{\longrightarrow}
Emb(\coprod_{q=1}^N B^4(\lambda_q);\cPtu,\sigma_{std})$$ is surjective. 
(where $Res(\coprod\varphi_{_q})=
\coprod(\varphi_{_q}\mid_{_{B^4(\lambda_q)}})$).\\
2) Let $1\leq r<N$, then 
${\dst Emb(\coprod_{q=1}^r B^4(\mu_q);\cPtu,\sigma_{std})}$ is 
connected  iff \\
$$Emb(\coprod_{q=1}^N B^4(\mu_q);\cPtu,\sigma_{std})\stackrel{Res}
{\longrightarrow}
Emb(\coprod_{q=1}^r B^4(\mu_q);\cPtu,\sigma_{std})$$ is surjective. 
(Where $Res(\coprod_1^N\varphi_{_q})=\coprod_1^r\varphi_{_q}$).

\end{prop}
 
\pf
1) Suppose that $Emb(\underline{\lambda})$ is connected 
($\underline{\lambda}$ stands for $(\lambda_1,\ldots,\lambda_N)$).\\
Let $\varphi_{_1}:\coprod_{q=1}^NB^4(\lambda_q)\rightarrow 
\CPtu\;\in Emb(\underline{\lambda})$.
Choose any $\bbar{\varphi}:\coprod_{q=1}^N(\mu_q)\rightarrow 
\CPtu\;\in Emb(\underline{\mu})$, 
and set $\varphi_{_0}=Res(\bbar{\varphi}) \in Emb(\underline{\lambda})$. 
According to our assumptions, $\varphi_{_0}$ and 
$\varphi_{_1}$ are isotopic through a symplectic isotopy 
$\{\varphi_{_t}\}_{t\in[0,1]}$. 
Let $L\in H_2(\CPtu;{\Bbb Z})$ be the homology class of a projective line.
Clearly there exists a smooth family of embedded 
2-spheres $S_t \subset (\CPtu,\sigstd)$ representing the homology class $L$,
such that for all $t$ \  $S_t\cap im\,\phit = \emptyset$. 
Using Banyaga's theorem on extensions of 
symplectic isotopies, extend $\varphi_{_t}$ to an ambient symplectic isotopy 
$F_t$. We have $F_t\circ\varphi_{_0}=\varphi_{_t}$ for all $t$, hence 
$F_1\circ\bbar{\varphi} \;\in Emb(\underline{\mu})$ and 
$Res(F_1\circ \bbar{\varphi})=\varphi_{_1}$.\\

Conversely, if $Res$ is surjective then clearly $Emb(\underline{\lambda})$ 
is connected.\\

\noindent 2) The proof is very similar to (1)'s.
\Qed

\nline
The proof of corollary \ref{cor-extend} follows now immediately by combining 
propositions \ref{prop-isotop-extend}~and~\ref{prop-pi0}.

\Qed

\section{Generalizations}

In this section we propose some generalizations of results stated
in section~\ref{sect-statement}.

\begin{prop}
\label{prop-pi0-gen1}

{\em (generalization of \ref{prop-pi0}(2))}

Let $0<\lambda_1\leq\lambda_2\leq\lambda_3\leq\lambda_4<\frac{1}{\sqrt{2}}$, 
and suppose that ${\dst \sum_{q=1}^4\lambda_q^2 < 1+2\lambda_1^2}$, then 
${\dst Emb(\coprod_{q=1}^4 B^4(\lambda_q); \cPtu,\sigma_{std})}$ is connected.

\end{prop}

\pf The proof is based on the inflation procedure, but here we have to 
overcome some technical difficulties arising from the balls being of 
different sizes.

As before, we would like to make inflation with curves in suitable homology 
classes to produce the forms $\tau^j_t$, and then find coefficients 
$x(t)>0\;y^j(t)\geq0$ satisfying 
\begin{equation}
x(t)(l-\sum_{q=1}^4\alpha_q(t)^2e_q) + \sum_{j=1}^r y^j(t)[\tau^j_t]
=l-\sum_{q=1}^4\lambda_q^2 e_q 
\label{eq:pi0-gen1-1}
\end{equation}
The problem is that in our case, if we attempt to make inflation with curves
in the class $A=2L-\sum_{q=1}^4E_q$ then equation (\ref{eq:pi0-gen1-1}) 
gives us a system of 5 equations in 2 unknowns $x(t),y(t)$, which generally 
does not have any solutions (unless $\lambda_1=\cdots=\lambda_4$ for example).

The idea is to control the values of the shrinking coefficients $\alpha_q(t)$ 
in such a way that the equations in this system become dependent.\\

Suppose that $\lambda_1\leq\cdots\leq\lambda_4<\frac{1}{\sqrt{2}}$ and 
$\sum_{q=1}^4\lambda_q^2<1+2\lambda_1^2$.
By making inflation with conics through the 4 blown-up points we get the forms 
$\tau_t$ in the cohomology class $[\tau_t]=2l-\sum_{q=1}^4e_q$. In our case
the cohomological condition (\ref{eq:pi0-gen1-1}) gives the system 
\begin{eqnarray}
x(t)+2y(t)=1 \label{eq:pi0-gen1-2} \\
\alpha_q(t)^2x(t)+y(t)=\lambda_q^2 \nonumber \\
q=1,\ldots,4 \nonumber
\end{eqnarray}

To make this system of 5 equations have a solution we must control the values 
of the $\alpha_q(t)$'s. It follows from the proof of the preparation
lemma that each of the $\alpha_q(t)$ can be chosen independently from the 
others, and that each $\alpha_q(t)$ can be chosen to be any positive smooth 
function satisfying the following conditions:\\
1) $0<\alpha_q(t)\leq\lambda_q$ for all t.\\
2) $\alpha_q(-1)=\alpha_q(1)=\lambda_q$.\\
3) $\alpha_q(0)$ is small enough (more precisely, we need that 
$\varphi_{_{-1}},\varphi_{_{1}}$ be holomorphic on $B^4(\alpha_q(0))$\ ).\\

In view of this, choose a small positive number $\epsilon$, and 
take $\alpha_1$ be a smooth function satisfying conditions 1-3 above, and such 
that $\alpha_1(t)\leq\epsilon$ for t near $0$, say $t\in[-\delta,\delta]$.
Now define $\alpha_2,\alpha_3,\alpha_4$ to be such that 
\begin{equation}
\label{eq:pi0-gen1-3}
\alpha_q(t)^2=
\frac{1}{2}[1-\frac{1-2\lambda_q^2}{1-2\lambda_1^2}(1-2\alpha_1(t)^2)]
\;\;\;q=2,3,4
\end{equation}

A straightforward computation shows that now the system 
(\ref{eq:pi0-gen1-2}) has solutions 
$x(t)>0,y^j(t)\geq0$, and that the $\alpha_q$'s satisfy conditions 1-2.
The only problem is that $\alpha_2,\alpha_3,\alpha_4$ may not 
satisfy the 3'd condition as one can check from equation (\ref{eq:pi0-gen1-3}).

To overcome this, note that when $\epsilon\rightarrow0^+$ we have 
$\alpha_q(\pm\delta)^2\rightarrow\frac{\lambda_q^2-\lambda_1^2}{1-2\lambda_1^2}
\;\;\;q=1,2,3,4$. Hence we have

$$\sum_{q=1}^4\alpha_q(\pm\delta)^2
{\raisebox{-1.72ex}
{$ \stackrel{\displaystyle \longrightarrow}
{\scriptstyle \epsilon\rightarrow 0^+}$}}
\frac{\sum_{q=1}^4\lambda_q^2\;\;-4\lambda_1^2}{1-2\lambda_1^2} < 1$$

Choose a small enough $\epsilon$, for which 
$\sum_{q=1}^4\alpha_q(\pm\delta)^2 < 1$.
>From proposition \ref{prop-N-balls} and it's proof it is easy to see that 
the forms $\Omgbar_{-\delta}$ and $\Omgbar_{\delta}$ are {\it {isotopic}}.

We have now a pseudo-isotopy connecting $\Omgbar_{-1}$ with 
$\Omgbar_1$, consisting of 3 parts:
$-1\leq t\leq-\delta$,\ \  $-\delta\leq t\leq\delta$,
\ \  $\delta\leq t\leq 1$, where the 2'nd part comes from the isotopy between 
$\Omgbar_{-\delta}$ and $\Omgbar_{\delta}$, while the 1'st and 3'd parts are 
taken from $\Omgbar_t$. Denote this pseudo-isotopy by $\Omgbar'_t$,
and suppose that $[\Omgbar'_t]=l-\sum_{q=1}^4 \alpha_q'(t)^2 e_q$.
According to the way we have constructed this pseudo-isotopy, if we replace 
the $\alpha_q$'s by $\alpha_q'$'s in the system (\ref{eq:pi0-gen1-2}), 
it becomes solvable, hence $\Omgbar_{-1}$ and $\Omgbar_1$ are isotopic. 
\Qed

\begin{prop}
\label{prop-pi0-gen2}
{\em (generalization of \ref{prop-pi0}(3))}

Let $0<\lambda_1\leq\cdots\leq\lambda_5$ and suppose that 
${\dst \sum_{q=1}^5\lambda_q^2 \geq 4\lambda_5^2}$, then \\ 
${\dst Emb(\coprod_{q=1}^5 B^4(\lambda_q); \cPtu,\sigma_{std})}$ is connected 
(provided it is non-empty).

\end{prop}

\pf Here we make inflation with curves in the classes
$$A_i:=2L-\sum_{q=1,q\neq i}^5 E_q, \;\;\;i=1,\ldots,5$$ 
to obtain closed 2-forms $\tau_t^{i}$ in the classes 
$$[\tau_t^{i}]=2l-\sum_{q=1,q\neq i}^5 e_q$$
The correction of the pseudo-isotopy to an isotopy yields the equations
\begin{eqnarray*}
x(t)+2\sum_{q=1}^5y^q(t)=1 \\
\alpha_{i-4}(t)^2 x(t)+y^{i-3}(t)+y^{i-2}(t)+y^{i-1}(t)+y^{i}(t)=
\lambda_{i-4}^2\\
i=1,\ldots,5
\end{eqnarray*}
(here we use the convention that $y^0(t)=y^5(t)\;\; y^{-1}=y^4(t)$ etc.).

This system has the solution
$$x(t)=\frac{2-\sum_{q=1}^5\lambda_q^2}{2-\sum_{q=1}^5\alpha_q(t)^2}$$
$$y^i(t)=\frac{1}{4}
[(\sum_{q=1}^5\lambda_q^2\;-4\lambda_i^2)-x(t)
(\sum_{q=1}^5\alpha_q(t)^2\;-4\alpha_i(t)^2)]$$

The proof is completed as soon as we show that $x(t)>0,\;y^i(t)\geq 0$ 
for all $t$.
Since $0<\alpha_q(t)\leq\lambda_q$ and $\sum_{q=1}^5\lambda_q^2 < 2$ 
(the last inequality comes from packing obstruction, see \cite{M-P}), 
we have $0<x(t)\leq 1$.
In order that the $y^i(t)$ be non-negative we need that 
\begin{equation}
\label{eq:pi0-gen2-1}
\sum_{q=1}^5\lambda_q^2\;-4\lambda_i^2\geq 
x(t)(\sum_{q=1}^5\alpha_q(t)^2\;-4\alpha_i(t)^2)
\end{equation}
for all $i,t$.

For this inequality to hold we control the values of the 
$\alpha_q$'s in the following way 
(compare with the proof of proposition \ref{prop-pi0-gen1}):

Choose some small positive number $\delta$.
When $\mid t\mid\geq\delta$ define $\alpha_q(t)=\lb \mid t\mid\lambda_q$, 
while on 
$\mid t\mid\leq\delta$ take the $\alpha_q(t)$'s to be such that the following
inequalities hold:\\
1) $0<\alpha_1(t)\leq\alpha_2(t)\leq\cdots\leq\alpha_5(t)$\\
2) $\sum_{q=1}^5\lambda_q^2\;-4\lambda_5^2\geq
\sum_{q=1}^5\alpha_q(t)^2\;-4\alpha_1(t)^2$\\
%
%
%
%

It is easily seen that such a controlling of the $\alpha_q$'s is possible 
when $\delta$ is small enough,   
and that the above conditions imply that inequality (\ref{eq:pi0-gen2-1}) 
holds, hence $y^i(t)\geq 0$.
\Qed

\nline
\noindent 
The following proposition can be proved in a similar way as 
\ref{prop-pi0-gen1} and \ref{prop-pi0}(1).

\begin{prop}
\label{prop-pi0-gen3}
{\em (generalization of \ref{prop-pi0}(1))}

Let $0<\lambda_1,\lambda_2,\lambda_3<\frac{1}{\sqrt{2}}$ then 
${\dst Emb(\coprod_{q=1}^3 B^4(\lambda_q); \cPtu,\sigma_{std})}$ is connected.

\end{prop}

\begin{prop}
\label{prop-pi0-gen4}
Let $N\geq 7,\; 0<\lambda<\sqrt{\frac{13}{5N}}$ then the space \\ 
${\dst Emb(\coprod_{1}^N B^4(\lambda); \cPtu,\sigma_{std})}$ is non-empty and 
connected.

\end{prop}
 
{\bf Outline of the proof \ \ } To prove this, one has to make inflation with 
$N$ quintics in $\CPtu_N$, in the homology classes 
$A_i=5L-2\sum_{q=i}^{i+5}E_q-E_{i+6}$ \ $i=1,\ldots,N$ 
(here we use the convention that 
$E_{N+r}=E_r$ for $r=1,\ldots,6$).
All the other details of the proof are 
quite similar to those of the proofs of propositions 
\ref{prop-pi0} and \ref{prop-N-balls}. Here again one has to overcome 
the non-positiveness of $\Jbar_s$. This can be done by using 
remark \ref{rem-holo-sphere}.

It remains to show that the above embedding spaces are non-empty.
This follows from a result due to Xu (see \cite{Xu}), according to which 
$Emb(\coprod_{1}^N B^4(\lambda); \lb \CPtu,\sigma_{std})\neq\emptyset$,
when $\lambda^2\leq\frac{\sqrt{N-1}}{N}$. 
Indeed, if $N\geq 8$, then $\frac{13}{5N} \leq\frac{\sqrt{N-1}}{N}$, 
while for $N=7$ it is well known (see \cite{M-P}) that there exists a 
symplectic embedding of 7 equal disjoint 4-balls of any radius
which is less then $\sqrt{\frac{3}{8}}$.
Clearly $\sqrt{\frac{13}{5\cdot 7}} < \sqrt{\frac{3}{8}}$. This proves 
that $Emb(\coprod_{1}^N B^4(\lambda); \CPtu,\sigma_{std})$ is non-empty 
when  $N\geq 7,\; 0<\lambda<\sqrt{\frac{13}{5N}}$.
\Qed

\begin{prop}
\label{prop-pi0-gen5}
Let $0<\lambda_1\leq\lambda_2\leq\lambda_3\leq\lambda_4<\lambda<
\frac{1}{\sqrt{2}}$, and suppose that ${\dst \sum_{q=1}^4\lambda_q^2<1+2\lambda_1^2}$,
then every symplectic embedding 
${\dst \coprod_{q=1}^4 B^4(\lambda_q)\rightarrow(\cPtu,\sigma_{std})}$ can be 
extended to a symplectic embedding 
${\dst \coprod_{1}^4  B^4(\lambda) \rightarrow(\cPtu,\sigma_{std})}$.

\end{prop}

\pf This is just a combination of proposition \ref{prop-pi0-gen1} and 
\ref{prop-isotop-extend}.
\Qed

\nline

By combining propositions \ref{prop-pi0-gen2},\ref{prop-pi0-gen3} with 
proposition \ref{prop-isotop-extend}, one can prove similar results for 
extensions of embeddings of 3 and of 5 balls of different radii.

\subsection{Remarks}

The inflation method can also be used to provide some information on higher 
homotopy groups of the embedding spaces. For example, in order to show that 
$\pi_1$ of some embedding space vanishes, one has to make 2-parametric 
inflation in order to span loops of symplectic embeddings by discs of such 
embeddings. However, one has to be careful when using 
pseudo-holomorphic curves, since generally in 2-parametric families 
bubbling with curves having 1'st Chern number zero may occur. 
We leave this discussion to another paper.

\nline
\subsection*{Acknowledgments}
I am deeply grateful to my teacher Prof. Leonid~Polterovich, who interested 
me in the whole subject, patiently spent many hours of explaining and 
listening, and made me aware of the ``inflation'' technique.

\end{document}